\documentclass[12pt]{article}
\usepackage{amsmath,amsthm,amssymb,color,epic,eepic,
cite}
\newlength{\minitwocolumn}
\newcommand{\Z}{{\Bbb Z}} 
\newcommand{\C}{{\Bbb C}} 
\newcommand{\FF}{{\Bbb F}} 
\newcommand{\K}{{\Bbb K}}

\newcommand{\cD}{{\cal D}}

\newcommand{\cR}{{\cal R}}

\newcommand{\cQ}{{\cal Q}}

\newcommand{\hL}{\widehat{L}}

\renewcommand{\H}{{\cal H}}
\newcommand{\la}{\lambda}

\newcommand{\al}{\alpha}

\newcommand{\vep}{\varepsilon}

\newcommand{\bH}{\bar{H}}

\newcommand{\hf}{\widehat{f}}

\newcommand{\noi}{{\noindent}}
\newcommand{\nn}{{\nonumber}}
\newcommand{\bea}{\begin{eqnarray}}
\newcommand{\ena}{\end{eqnarray}}
\newcommand{\beit}{\begin{itemize}}
\newcommand{\enit}{\end{itemize}}

\newcommand{\be}{\begin{eqnarray*}}
\newcommand{\en}{\end{eqnarray*}}
\newcommand{\lb}[1]{\label{#1}}

\newcommand{\ds}[1]{{\displaystyle #1 }}

\newcommand{\End}{{\rm End}}

\newcommand{\id}{{\rm id}}

\def\infq4p#1{{(#1;q^4,p)_\infty}}

\newcommand{\hPhi}{\widehat{\Phi}}

\newcommand{\hPsi}{\widehat{\Psi}}

\newcommand{\hLp}{\widehat{L}^+}

\newcommand{\tot}{\widetilde{\otimes}}

\newcommand{\mmatrix}[1]{\begin{matrix} #1 \end{matrix}}

\font\teneufm=eufm10
\font\seveneufm=eufm7
\font\fiveeufm=eufm5
\newfam\eufmfam
\textfont\eufmfam=\teneufm
\scriptfont\eufmfam=\seveneufm
\scriptscriptfont\eufmfam=\fiveeufm

\let\goth\frak
\newcommand{\slth}{\widehat{\goth{sl}}_2}
\newcommand{\slt}{\goth{sl}_2}

\newcommand{\Bqla}{{{\cal B}_{q,\lambda}}}
\newcommand{\Uqp}{U_{q,p}}

\newcommand{\hh}{\goth{h}}
\newcommand{\h}{H}

\font\seventeeneufm=eufm10 scaled\magstep3   
   
\newcommand{\slthBig}{\widehat{\mbox{\seventeeneufm sl}}_2} 
\makeatletter
\@addtoreset{equation}{section}
\makeatother

\newtheorem{thm}{Theorem}[section]
\newtheorem{prop}[thm]{Proposition}
\newtheorem{lem}[thm]{Lemma}

\newtheorem{dfn}[thm]{Definition}

\begin{document}
\bibliographystyle{unsrt}

\vspace{2cm}
\begin{center}
{\Large \bf Elliptic Quantum Group $U_{q,p}(\slthBig)$\\
  and Vertex Operators\\[10mm] }
{\large  Hitoshi Konno}\\[6mm]
{\it Department of Mathematics, Graduate School of Science, 
\\Hiroshima University, Higashi-Hiroshima 739-8521, Japan\\
       konno@mis.hiroshima-u.ac.jp}\\[10mm]
\end{center}

\begin{abstract}
\noindent 
Introducing an $\h$-Hopf algebroid structure into $U_{q,p}(\slth)$, we
 investigate the vertex operators of the elliptic quantum group 
$U_{q,p}(\slth)$ defined as  
intertwining operators of infinite dimensional 
$U_{q,p}(\slth)$-modules. We show that the vertex operators coincide with 
the previous results obtained indirectly 
by using the quasi-Hopf algebra $\Bqla(\slth)$. 
This shows a consistency of our $\h$-Hopf algebroid structure
 even in the case with non-zero central element. 

\end{abstract}
\nopagebreak

\section{The Elliptic Algebra $U_{q,p}(\slthBig)$}
In this section we review a definition of the elliptic algebra $U_{q,p}(\slth)$
and its $RLL$ formulation following \cite{Konno,JKOS2}.
\subsection{Definition of $U_{q,p}(\slthBig)$}\lb{app:1.3}
The elliptic algebra 
$\Uqp(\slth)$ was introduced in \cite{Konno} as an elliptic analogue of 
the quantum affine algebra $U_q(\slth)$ in the Drinfeld realization. 
It was soon realized that $\Uqp(\slth)$ is isomorphic to the  
 the tensor product of $U_q(\slth)$  and a Heisenberg algebra 
$\{P,e^Q\}$\cite{JKOS2}. We here define $\Uqp(\slth)$ along
 the latter observation. 

Let us fix a complex number $q$ such that $q\not=0, |q|<1$. 
\begin{dfn}\cite{Drinfeld}\lb{defUq}
For a field $\K$, the quantum affine algebra $\K[U_q(\slth)]$ 
in the Drinfeld realization  is an associative algebra over $\K$ 
generated by the Drinfeld generators 
$a_n\ (n\in \Z_{\not=0})$, $x_n^\pm\ (n\in \Z)$, 
$h$,  $c, d$.    
The defining relations are given as follows.
\be
&&c :\hbox{ central },\nn\\
&& [h,d]=0,\quad [d,a_{n}]=n a_{n},\quad 
[d,x^{\pm}_{n}]=n x^{\pm}_{n}, \nn\\
&&[h,a_{n}]=0,\qquad [h, x^\pm(z)]=\pm 2 x^{\pm}(z),\nn\\
&&
[a_{n},a_{m}]=\frac{[2n]_{q}[c n]_{q}}{n}
q^{-c|n|}\delta_{n+m,0},\nn\\
&&
[a_{n},x^+(z)]=\frac{[2n]_{q}}{n}q^{-c|n|}z^n x^+(z),\nn\\
&&
[a_{n},x^-(z)]=-\frac{[2n]_{q}}{n} z^n x^-(z),
\lb{uqslth}\\
&&(z-q^{\pm 2}w)
x^\pm(z)x^\pm(w)= (q^{\pm 2}z-w) x^\pm(w)x^\pm(z),\nn
\\
&&[x^+(z),x^-(w)]=\frac{1}{q-q^{-1}}
\left(\delta\bigl(q^{-c}\frac{z}{w}\bigr)\psi(q^{c/2}w)
-\delta\bigl(q^{c}\frac{z}{w}\bigr)\varphi(q^{-c/2}w)
\right).\nn
\en
where $[n]_q=\frac{q^n-q^{-n}}{q-q^{-1}}$, $\delta(z)=\sum_{n\in\Z}z^n$ and 
\be
&&x^\pm(z)=\sum_{n\in \Z}x^\pm_{n} z^{-n},\\
&&\psi(q^{c/2}z)=q^{h}
\exp\left( (q-q^{-1}) \sum_{n>0} a_{n}z^{- n}\right),\quad 
\varphi(q^{-c/2}z)=q^{-h}
\exp\left(-(q-q^{-1})\sum_{n>0} a_{-n}z^{ n}\right).
\en
\end{dfn}

Let $r$ be a complex parameter. We set $r^*=r-c$, $p=q^{2r}$ and $p^*=q^{2r^*}$. 
We define the  Jacobi theta functions $[u]$ and $[u]^*$ by
\be
&&[u]=\frac{q^{{u^2}/{r}-u}}{(p;p)_\infty^3}\Theta_{p}(q^{2u}),\qquad 
[u]^*=\frac{q^{{u^2}/{r^*}-u}}{(p^*;p^*)_\infty^3}\Theta_{p^*}(q^{2u}),
\en
where
\be
&&\Theta_{p}(z)=(z;p)_\infty(p/z;p)_\infty(p;p)_\infty,\\
&&(z;p_1,p_2,\cdots, p_m)_\infty=\prod_{n_1,n_2,\cdots,n_m=0}^\infty(1-zp_1^{n_1}p_2^{n_2}\cdots p_m^{n_m}).
\en
Setting $p=e^{-{2\pi i}/{ \tau}}$, $[u]$ satisfies the quasi-periodicity $[u+r]=-[u]$, $[u+r\tau]=e^{-\pi i (2u/r+\tau)}[u]$.

We denote by $\{P,e^Q\}$ a Heisenberg algebra 
commuting with $\C[U_q(\slth)]$ and satisfying  
\bea
&&[P, e^{Q}]=-e^{Q}. \lb{Heisenberg}
\ena
We take the realization $Q=\frac{\partial }{\partial P}$. 
We set $\h=\C P\oplus \C  r^*$ 
 and $\h^*= \C Q\oplus \C \frac{\partial}{\partial r^*}$
 with the  pairing $<\ ,\ >$
 \be
&&<Q,P>=1=<\frac{\partial}{\partial r^*},r^*>,
\en
the others are zero. 

We also consider the Abelian group $\bar{\h}^{*}= \Z Q$. 
We denote by $\C[\bar{H}^*]$ the group algebra over $\C$ of $\bar{H}^*$,
 and by $e^{\al}$ the element of $\C[\bar{H}^*]$ corresponding to $\al\in \bar{H}^*$. 
These $e^\al$ satisfy $e^\al e^\beta=e^{\al+\beta}$
 and $(e^\al)^{-1}=e^{-\al}$. 
In particular, $e^0=1$ is the identity element. 

Now we take the power series field  $\FF=\C((P,r^*))$ 
as $\K$ and consider the semi-direct product $\C$-algebra 
$U_{q,p}(\slth)=\FF[U_q(\slth)]\otimes_{\C} \C[\bar{H}^*]$ of 
$\FF[U_q(\slth)]$ and 
$\C[H^*]$, whose multiplication is defined by
\be
&&(f(P,r^*)a\otimes e^{\al})\cdot(g(P,r^*)b\otimes e^{\beta})=f(P,r^*)g(P+<\al,P>,r^*)ab\otimes e^{\al+\beta},\\
&&a, b\in \C[U_q(\slth)],\  f(P,r^*), g(P,r^*)\in \FF,\ \al, \beta\in \bH^*. 
\en

Let us consider the following generating functions. 
\be
&&
u^+(z,p)=\exp\left(\sum_{n>0}\frac{1}{[r^*n]_q}a_{-n}(q^rz)^n\right),\quad 
u^-(z,p)=\exp\left(-\sum_{n>0}\frac{1}{[rn]_q}a_{n}(q^{-r}z)^{-n}\right).
\en 
We define an automorphism $\phi_r$ of $\C[U_q(\slth)]$ by
\be
&&c\mapsto c,\quad h\mapsto h, \quad d\mapsto d,\\
&&x^+(z)\mapsto  u^+(z,p)x^+(z),\quad x^-(z)\mapsto x^-(z)u^-(z,p),\\
&&\psi(z)\mapsto u^+(q^{c/2}z,p)\psi(z)u^-(q^{-c/2}z,p),\\
&&\varphi(z)\mapsto u^+(q^{-c/2}z,p)\varphi(z)u^-(q^{c/2}z,p).
\en

\begin{dfn}\lb{ecurrents}
We define  $E(u), F(u), K(u)\in U_{q,p}(\slth)[[u]]$ and $\hat{d}$ by 
the following formulae.
\be
\!\!\! E(u)&=&\phi_r(x^+(z))e^{2Q}z^{-{(P-1)}/{r^*}},\\
\!\!\! F(u)&=&\phi_r(x^-(z))z ^{{(P+ h-1)}/{r}},\\
\!\!\!K(u)&=&\exp\left(\sum_{n>0}\frac{[n]_q}{[2n]_q[r^*n]_q}a_{-n}(q^cz)^n\right)
\exp\left(-\sum_{n>0}\frac{[n]_q}{[2n]_q[rn]_q}a_{n}z^{-n}\right)\\
&&\qquad\qquad\times e^Qz^{-{c}(2P-1)/{4rr^*}+h/{2r}},
\\
\!\!\!\hat{d}&=&d-\frac{1}{4r^*}(P-1)(P+1)+\frac{1}{4r}(P+h-1)(P+h+1), 
\en
where we set $z=q^{2u}$. 
We call $E(u), F(u)$,  $K(u)$ the elliptic currents.
\end{dfn}
In fact, from Definition \ref{defUq} and \eqref{Heisenberg}, we can derive the following relations.
\begin{prop}\lb{Defrels}
\be
&&c:\hbox{ {\rm central}}, 
\lb{u1}
\\
&&[h,a_{n}]=0,\quad [h,E(u)]=2E(u),\quad 
[h,F(u)]=-2 F(u), 
\lb{u2}\\
&&[\hat{d},h]=0,\quad [\hat{d},a_{n}]= n a_{n},\quad\\ 
&&[\hat{d},E(u)]=\left(-z\frac{\partial}{\partial z}-\frac{1}{r^*}\right)
E(u) , \quad
[\hat{d},F(u)]=\left(-z\frac{\partial}{\partial z}-\frac{1}{r}\right)
F(u) ,
\lb{u3}\\
&&
[a_{n},a_{m}]=\frac{[2n]_q[c n]_q}{n}
q^{-c|n|}\delta_{n+m,0},
\label{u4}\\
&&
[a_{n},E(u)]=\frac{[2n]_q}{n}q^{-c|n|}z^n E(u),
\lb{u5}\\
&&
[a_{n},F(u)]=-\frac{[2n]_q}{n} z^n F(u),
\lb{u6}\\
&& E(u)E(v)
=\frac{\left[u-v+1\right]^*}{\left[u-v-1\right]^* } 
E(v)E(u),
\lb{u7}\\
&&  F(u)F(v)
=\frac{\left[u-v-1\right]}{\left[u-v+1 \right]}
F(v)F(u),
\lb{u8}\\
&&[E(u),F(v)]
=\frac{1}{q-q^{-1}}
\left(\delta\left(q^{-c}\frac{z}{w}\right)H^+(q^{c/2}w)
-\delta\left(q^{c}\frac{z}{w}\right)H^-(q^{-c/2}w)
\right),\lb{u9}
\en
where $z=q^{2u}$, $w=q^{2v}$, 
\be
&&H^\pm(z)=\kappa K\left(u\pm\frac{1}{2}(r-\frac{c}{2})+\frac{1}{2}\right)
 K\left(u\pm\frac{1}{2}(r-\frac{c}{2})-\frac{1}{2}\right),\lb{HKK}\\
&&\kappa=\lim_{z\to q^{-2}}\frac{\xi(z;p^*,q)}{\xi(z;p,q)},\qquad 
\xi(z;p,q)=\frac{(q^2z;p,q^4)_\infty(pq^2z;p,q^4)_\infty}
{(q^4z;p,q^4)_\infty(pz;p,q^4)_\infty}.\nn
\en
\end{prop}
\noi
In particular we have the following relations which, together with the last three relations in the above, appeared in \cite{Konno}. 
\begin{prop}
\be
\!\!\!\!\!\!\!\!\!\!&&K(u)K(v)=\rho(u-v)K(v)K(u),
\lb{Ucom1}\\[2mm]
\!\!\!\!\!\!\!\!\!\!&&K(u)E(v)=\frac{[u-v+\frac{1-r^*}{2}]^*}
{[u-v-\frac{1+r^*}{2}]^*}E(v)K(u),
\lb{Ucom2}\\[2mm]
\!\!\!\!\!\!\!\!\!\!&&K(u)F(v)=\frac{[u-v-\frac{1+r}{2}]}
{[u-v+\frac{1-r}{2}]}F(v)K(u),
\lb{Ucom3}\\
&&H^{+}(u)H^{-}(v)=\frac{[{u-v-1-\frac{c}{2}}]}{
[{u-v+1-\frac{c}{2}}]}
\frac{[{u-v+1+\frac{c}{2}}]^*}{
[{u-v-1+\frac{c}{2}}]^*}H^{-}(v)H^{+}(u),\\
&&H^{\pm}(u)H^{\pm}(v)=\frac{[{u-v-1}]}{
[{u-v+1}]}
\frac{[{u-v+1}]^*}{
[{u-v-1}]^*}H^{\pm}(v)H^{\pm}(u),
\en
where
\be
&&\rho(u)=\frac{\rho^{+*}(u)}{\rho^+(u)},\quad
\rho^+(u)=z^{{1}/{2r}}
\frac{\{pq^2z\}^2}{\{pz\}\{pq^4z\}}
\frac{\{z^{-1}\}\{q^4z^{-1}\}}{\{q^2z^{-1}\}^2},\quad \{z\}=(z;p,q^4)_\infty,\\
&&\rho^{+*}(u)=\rho^+(u)|_{r\to r^*}.
\en
\end{prop}

\begin{dfn}
We call a set ( $\FF[U_q(\slth)]\otimes_{\C}\C[\bH^*]$,  $\phi_r$)  
the elliptic algebra $U_{q,p}(\slth)$.
\end{dfn}
The following relations are also useful. 
\begin{prop}\lb{commPh}
\be
&&[K(u),P]=K(u),\quad [E(u),P]=2E(u),\quad [F(u),P]=0,
\\
&&[K(u),P+h]=K(u),\quad [E(u),P+h]=0,\quad [F(u),P+h]=2F(u).
\en
\end{prop}

\subsection{The $RLL$-relation for $U_{q,p}(\slthBig)$}

We next summarize the $RLL$-relation for $U_{q,p}(\slth)$\cite{JKOS2}. 
 In the next section, the $L$ operator is used to 
discuss the $\h$-Hopf algebroid structure of $U_{q,p}(\slth)$.

Let us define the half currents in the following way.
\begin{dfn}\lb{halfcurrents}
\be
&&
K^+(u)=K(u+\tfrac{r+1}{2}),
\label{kplush}\\
&&E^+(u)
=a^* \oint_{C^*} E(u') 
\frac{\left[u-u'+c/2-P+1\right]^*[1]^*}
{[u-u'+c/2]^*[P-1]^*}
\frac{dz'}{2\pi i z'},
\lb{Eplus}\\
&&F^+(u)
=a \oint_{C} F(u') 
\frac{\left[u-u'+P+h-1\right][1]}{[u-u'][P+h-1]}
\frac{dz'}{2\pi i z'}. 
\lb{Fplus}
\en
Here the contours are chosen such that  
\be
&&C^* : |p^*q^c z|<|z'|<|q^cz|, \qquad
C : |pz|<|z'|<|z|,
\lb{C}
\en
and the constants $a,a^*$ are chosen to satisfy 
$
{a^* a [1]^*\kappa\over q-q^{-1}} =1.
$
\end{dfn}

\begin{dfn}\lb{Loperator}
We define the operator $\hL^+(u)\in \End_{\C} V\otimes\Uqp\bigl(\slth\bigr)$ 
with $V\cong\C^2$, 
by 
\be
&&\hL^+(u)=
\left(
\begin{array}{cc}
1 &F^+(u) \\
0 &1            \\
\end{array}
\right)
\left(
\begin{array}{cc}
K^+(u-1) & 0 \\
0             &K^+(u)^{-1}\\
\end{array}
\right)
\left(
\begin{array}{cc}
1            &0      \\
E^+(u) &1      \\
\end{array}
\right).
\lb{Gauss}
\en
\end{dfn}

\begin{prop}\lb{prop:RLL1}
The  operator $\hL^+(u)$ satisfies the following $RLL$ relation. 
\bea
&&\hspace*{-15mm}R^{+(12)}(u_1-u_2,P+h)
\hL^{+(1)}(u_1)\hL^{+(2)}(u_2)
=
\hL^{+(2)}(u_2)\hL^{+(1)}(u_1)
R^{+*(12)}(u_1-u_2,P),
\lb{RLLrel}
\ena
where $R^{+}(u,P+h)$ and $R^{+*}(u,P)=R^{+}(u,P)|_{r\to r^*}$ denote the elliptic 
dynamical $R$ matrices given by 
\bea
&&R^+(u,s)=\rho^+(u)
\left(
\begin{array}{cccc}
1 &                  &             & \\
  &b(u,s)      &c(u,s) & \\
  &\bar{c}(u,s)&\bar{b}(u,s) & \\
  &                  &             &1 \\
\end{array}
\right)
\lb{Rmat}
\ena
with 
\be
&&b(u,s)=
\frac{[s+1] [s-1]  }{[s]^2}
\frac{[u]}{[1+u]},
\quad 
c(u,s)=
\frac{[1]}{[s]}
\frac{[s+u]}{[1+u]},
\lb{Rmat1}\\
&&\bar{c}(u,s)=\frac{[1]}{[s]}
\frac{[s-u]}{[1+u]},
\quad \qquad\qquad\quad\; 
\bar{b}(u,s)=
\frac{[u]}{[1+u]}.
\lb{Rmat4}
\en
\end{prop}

Note that if we set $L^+(u,P)=\hL^+(u)e^{-h\otimes Q}$, $L^+(u,P)$ is 
independent of $Q$ and satisfies the dynamical $RLL$ relation\cite{JKOS2} 
characterizing the quasi-Hopf algebra $\Bqla(\slth)$\cite{JKOS}. 
Moreover with the parametrization $\la=(r^*+2)\Lambda_0+(P+1)\bar{\Lambda}_1$, where $\Lambda_0, \Lambda_0+\bar{\Lambda}_1$ are the fundamental weights of $\slth$,  $\Bqla(\slth)$ is isomorphic to $\FF[U_q(\slth)]$, as an associative algebra. These two facts lead to the isomorphism $U_{q,p}(\slth)\cong \Bqla(\slth)\otimes_{\C}\C[\bH^*]$ as a semi-direct product $\C$-algebra. However this semi-direct product breaks down the quasi-Hopf algebra structure, so that $U_{q,p}(\slth)$ is not a quasi-Hopf algebra. In the next section, we show that a relevant coalgebra structure of $U_{q,p}(\slth)$ is 
the $H$-Hopf algebroid. 

 Note also that the $c=0$ case of the 
dynamical $RLL$ relation for $\Bqla(\slth)$ coincides with the 
one studied by Felder \cite{Felder,EF}, whereas the $c=0$ case of 
\eqref{RLLrel} coincides with  
the $RLL$ relation studied in 
\cite{EV1,EV2,KR} for the trigonometric $R$ and 
in \cite{KNR} for the elliptic $R$.

\section{$H$-Hopf Algebroid Structure of $U_{q,p}(\slthBig)$}
In this section we introduce an $\h$-Hopf algebroid structure into  
$U_{q,p}(\slth)$. The detailed discussion will 
be published elsewhere\cite{Konno07}. We follow the definition of 
$\h$-Hopf algebroid given in \cite{EV1,EV2} and \cite{KR,KNR} with 
a modification which makes it applicable in
 the case with non-zero central element. 

Let $\bar{\hh}=\C h$ be the Cartan subalgebra, $\al_1$ the simple root and 
$\bar{\Lambda}_1$ the fundamental weight of $\slt$. We set $\cQ=\Z \al_1$ and $\bar{\hh}^*=\C \bar{\Lambda}_1$. Let us use the same symbol $<,>$ to denote the standard paring of $\bar{\hh}$ and 
$\bar{\hh}^*$.  
Using the isomorphism  $\phi:\cQ\to \bar{\h}^*$ by $n\al_1\mapsto nQ$,  
we define 
the $\bar{\h}^*$-bigrading structure of $U_{q,p}=U_{q,p}(\slth)$ by 
\bea
&&{U}_{q,p}=\bigoplus_{\al, \beta \in \bar{\h}^*} (U_{q,p})_{\al \beta},\nn\\
&&(U_{q,p})_{\al\beta }=\left\{\ x\in U_{q,p}\ \left|\ \mmatrix{q^{h}xq^{-h}=q^{<\bar{\al},h>}x,\ \al=\phi(\bar{\al})+
\beta\cr
q^{P}xq^{-P}=q^{<\beta,P>}x \cr}\ \right.\right\}.\lb{bigrading}
\ena
Noting $<\bar{\al},h>=<\phi(\bar{\al}),P>$, we have $q^{P+h}xq^{-(P+h)}=q^{<\al,P>}x$ for $x\in (U_{q,p})_{\al\beta }$. 

We regard 
$\widehat{f}=f(P,r^*)\in \FF$ as a meromorphic function on $\h^*$ by 
\be
&&\widehat{f}(\mu)=f(<\mu,P>,<\mu,r^*>)\quad \mu\in {H}^*
\en
and consider  the field of  meromorphic functions  $M_{{\h}^*}$ on ${\h}^*$ given by 
\be
&&{M}_{{\h}^*}=\left\{ \widehat{f}:{\h}^*\to \C\ \left|\ 
\widehat{f}=f(P,r^*)\in \FF\right.\right\}.
\en 
We define two embeddings (the left and right moment maps) 
$\mu_l, \mu_r : {M}_{{\h}^*}\to (U_{q,p})_{00}$  by
\bea
\mu_l(\widehat{f})=f(P+h,{r^*+c}),\qquad 
\mu_r(\widehat{f})=f(P,{r^*}).\lb{mmUqp}
\ena
From \eqref{bigrading}, one finds for $x\in (U_{q,p})_{\al\beta}$
\be
&&\mu_l(\widehat{f})x=f(P+h,r^*+c)x=x f(P+h+<\al,P>,r^*+c)=x \mu_l({T}_{\al}
\widehat{f}),\lb{mlUqp}\\
&&\mu_r(\widehat{f})x=f(P,r^*)x=x f(P+<\beta,P>,r^*)=x \mu_r({T}_{\beta}\widehat{f}),\qquad \lb{mrUqp}
\en
where we regard ${T}_{\al}=e^\al\in \C[\bH^*]$ as 
a shift operator $M_{{\h}^*}\to M_{{\h}^*}$ 
\be
({T}_{\al}\widehat{f})=e^{\al}f(P,r^*)e^{-\al}={f}(P+<\al,P>,r^*).
\en
Hereafter we abbreviate 
$f(P+h,{r^*+c})$ and $f(P,{r^*})$ as $f(P+h)$ and
 $f^*(P)$, respectively.  

Then equipped with the 
bigrading structure \eqref{bigrading} and two moment 
maps \eqref{mmUqp}, the elliptic algebra $U_{q,p}(\slth)$ 
is an $\h$-algebra\cite{EV1,EV2}.

In addition, we need the $\h$-algebra $\cD$ of the shift operators given by
\be
&&\cD=\{\ \sum_i \widehat{f}_i{T}_{\al_i}\ |\ \widehat{f}_i \in M_{{\h}^*}, \al_i\in \bH^*\ \},\\
&&(\cD)_{ \al \al}=\{\ \widehat{f}{T}_{-\al}\ \},\quad 
(\cD)_{\al \beta}=0\quad \al\not=\beta,\\
&&\mu_l^{\cD}(\widehat{f})=\mu_r^{\cD}(\widehat{f})=\widehat{f}{T}_0 \qquad \widehat{f}\in M_{{\h}^*}.
\en

Let $A$ and $B$ be two $\h$-algebras, $U_{q,p}\ {\rm or}\ \cD$. 
The tensor product $A {\widetilde{\otimes}}B$ is the bigraded vector space with 
\be
 (A {\widetilde{\otimes}}B)_{\al\beta}=\bigoplus_{\gamma\in \bH^*} (A_{\al\gamma}\otimes_{M_{{\h}^*}}B_{\gamma\beta}),
\en
where $\otimes_{M_{{\h}^*}}$ denotes the usual tensor product 
modulo the following 
relations.
\bea
\mu^A_r(\widehat{f})a\otimes b=a\otimes \mu_l^B(\widehat{f})b\qquad a \in A, b \in B.\lb{deftot}
\ena
Then the tensor product $A {\widetilde{\otimes}}B$ is again an $\h$-algebra with the multiplication $(a\otimes b)(c\otimes d)=ac\otimes bd$ and the moment maps 
\be
\mu_l^{A {\widetilde{\otimes}}B} =\mu_l^A\otimes 1,\qquad \mu_r^{A {\widetilde{\otimes}}B} =1\otimes \mu_r^B.
\en
Note that we have the $\h$-algebra isomorphism $U_{q,p}\tot\cD\cong U_{q,p}\cong \cD\tot U_{q,p}$ by $x\tot T_{-\beta}=x=T_{-\al}\tot x$ for 
$x\in (U_{q,p})_{\al\beta}$. 

Now let us define an $\h$-Hopf algebroid structure on $U_{q,p}$  
as its co-algebra structure. 
For this purpose, it is convenient to 
use the $L$ operator $\hL^+(u)$. 
We shall write the entries of $\widehat{L}^+(u)$ as 
\be
&&\hLp(u)=\left(\mmatrix{\hLp_{++}(u)&\hLp_{+-}(u)\cr
\hLp_{-+}(u)&\hLp_{--}(u)\cr
}\right).
\en
From Proposition \ref{commPh} and Definition \ref{Loperator},  
 one finds 
\be
&&\hL^+_{\vep_1\vep_2}(u)\in (U_{q,p})_{-\vep_1Q,-\vep_2Q}.
\lb{grL1}
\en
It is also easy to check the relations  
\be
f(P+h)\hL^+_{\vep_1\vep_2}(u)&=&
\hL^+_{\vep_1\vep_2}(u)f(P+h-\vep_1),\\
f^*(P)\hL^+_{\vep_1\vep_2}(u)&=&\hL^+_{\vep_1\vep_2}(u)
f^*(P-\vep_2).
\en

\begin{dfn}\lb{Hopfalgebroid}
We define $\h$-algebra homomorphisms, $\vep : U_{q,p}\to \cD$ and  $\Delta : U_{q,p}\to U_{q,p} \widetilde{\otimes}U_{q,p}$ by
\be
&&\vep(\hL^+_{\vep_1\vep_2}(u))=\delta_{\vep_1,\vep_2}{T}_{-\vep_2 Q},\quad \vep(e^Q)=e^Q,
\lb{counitUqp}\\
&&\vep(\mu_l(\widehat{f}))=\vep(\mu_r(\widehat{f}))=\widehat{f}T_0, \lb{counitf}\\
&&\Delta(\hL^+_{\vep_1\vep_2}(u))=\sum_{\vep'}\hL^+_{\vep_1\vep'}(u)\widetilde{\otimes}
\hL^+_{\vep'\vep_2}(u),\lb{coproUqp}\\
&&\Delta(e^{Q})=e^{Q}\tot e^{Q},\\
&&\Delta(\mu_l(\hf))=\mu_l(\hf)\widetilde{\otimes} 1,\quad 
\Delta(\mu_r(\hf))=1\widetilde{\otimes} \mu_r(\hf).\lb{coprof}
\en
We also define an $\h$-algebra anti-homomorphism  $S : U_{q,p}\to U_{q,p}$ by
\be
&&S(\hL^+_{++})=\hL^+_{--}(u-1), \quad S(\hL^+_{+-}(u))=-\frac{[P+h+1]}{[P+h]}\hL^+_{+-}(u-1), \quad \\
&&S(\hL^+_{-+}(u))=-\frac{[P]^*}{[P+1]^*}\hL^+_{-+}(u-1),\quad S(\hL^+_{--}(u))=\frac{[P+h+1][P]^*}{[P+h][P+1]^*}\hL^+_{++}(u-1),\\
&&S(e^Q)=e^{-Q},\qquad S(\mu_r(\hf))=\mu_l(\hf),\quad S(\mu_l(\hf))=\mu_r(\hf).
\en
\end{dfn}

In fact one can show that $\Delta$ and $S$ preserve the $RLL$ 
relation \eqref{RLLrel}. 
Moreover we have the following lemma indicating that 
$\vep, \Delta$ and $S$ satisfy the axioms for 
the counit, the comultiplication and the antipode. 
Hence 
the $\h$-algebra $U_{q,p}(\slth)$ with $(\Delta,\vep,S)$ is an 
$\h$-Hopf algebroid\cite{EV1,EV2,KR}.
\begin{lem}\lb{counitcopro}
The maps $\vep$, $\Delta$ and $S$ satisfy
\be
&&(\Delta\otimes \id)\circ \Delta=(\id \otimes \Delta)\circ \Delta,\lb{coaso}\\
&&(\vep \otimes \id)\circ\Delta =\id =(\id \otimes \vep)\circ \Delta.\lb{vepDelta}\\
&&m\circ (\id \otimes S)\circ\Delta(x)=\mu_l(\vep(x)1),\quad \forall x\in U_{q,p},\\
&&m\circ (S\otimes\id  )\circ\Delta(x)=\mu_r(T_{\al}(\vep(x)1)),\quad \forall x\in (U_{q,p})_{\al \beta}.
\en
\end{lem}
\noi
\begin{dfn}
We call the $\h$-Hopf algebroid $(U_{q,p}(\slth),\h,M_{\h^*},\mu_l,\mu_r,\Delta,\vep,S)$ the  elliptic quantum group $U_{q,p}(\slth)$. 
\end{dfn}


\section{Representations}
We consider the dynamical 
representations, i.e. the representations as $\h$-algebras\cite{FV,EV1,EV2}, 
of the elliptic algebra $U_{q.p}(\slth)$.

\subsection{Evaluation Representation}
We construct the evaluation representation of $U_{q,p}(\slth)$  by 
using the one  of $\FF[U_q(\slth)]$.  
We define the $l+1$-dimensional vector space over $\FF$ by 
$\ds{V^{(l)}=\bigoplus_{m=0}^l \FF v^l_m}$. Here  $v^l_m\ (0\leq m\leq l)$ denote the 
weight vectors satisfying $h v^{l}_m =(l-2m)v^l_m$. 
Consider the operator $S^\pm$ acting on $V^{(l)}$ by
$S^{\pm}v^l_m=v^l_{m\mp1}, \ v^l_m=0\quad {\rm for } \ m<0,\ \ m>l.$
In terms of the Drinfeld generators, the evaluation representation $(\pi_{l,w},
 V^{(l)}_w=V^{(l)}\otimes \C[w,w^{-1}])$ of $\FF[U_q(\slth)]$ is given by\cite{JKOS2}
\be
&&\pi_{l,w}(c)=0,\qquad \pi_{l,w}(d)=0,\\
&&\pi_{l,w}(a_n)=\frac{w^n}{n}\frac{1}{q-q^{-1}}((q^n+q^{-n})q^{nh}-(q^{(l+1)n}+q^{-(l+1)n})),\\
&&\pi_{l,w}(x^{\pm}(z))=S^{\pm}\left[\frac{\pm h+l+2}{2}\right]_{q}\delta
\left(q^{h\pm1}\frac{w}{z}\right).
\en 
Note that  $\ds{V^{(l)}_w=\bigoplus_{\mu\in \{-l,-l+2,\cdots,l\}}V_\mu}$ 
with $V_\mu,\ \mu=l-2m$ spaned by $v^l_m\otimes w^n\ (n\in \Z)$. 

Let us define the $\h$-algebra $\cD_{\h,V}$ by
\be
&&\cD_{\h,V}=\bigoplus_{\al,\beta\in \bH^*}(\cD_{\h,V})_{\al\beta},\\
&&\hspace*{-10mm}(\cD_{\h,V})_{\al\beta}=
\left\{\ X\in \End_{\C}V^{(l)}_w\ \left|\ \mmatrix{X(f^*(P)v)=f^*(P-<\beta,P>)X(v), v\in V^{(l)}_w\cr 
X(V_\mu)\subseteq V_{\mu+\phi^{-1}(\al)-\phi^{-1}(\beta)},\ f^*(P)\in \FF\cr}\ \right.\right\},\\
&&\mu_l^{\cD_{H,V}}(\widehat{f})v=f(P+\mu)v,\quad 
\mu_r^{\cD_{H,V}}(\widehat{f})v=f^*(P)v\qquad 
\en
for $v\in V_{\mu}$. Then $\widehat{\pi}_{l,w}=\pi_{l,w}\otimes \id : 
U_{q,p}(\slth)=\FF[U_q(\slth)]\otimes_{\C} \C[\bH^*] \to \cD_{\h,V}$ with 
$e^Qv^l_m=v^l_m$ yields the $\h$-algebra homomorphism. 
We call $(\widehat{\pi}_{l,w},V^{(l)}_w)$ the dynamical 
evaluation representation. 
In particular, applying this to Definitions \ref{ecurrents}, \ref{halfcurrents}, \ref{Loperator}, we obtain the following 
expressions for the images of the $\hL^+(u)$ operator.
\begin{thm}\lb{repL}
\be
\widehat{\pi}_{l,w}(\hL^+_{++}(u))&=&-\frac{[u-v+\frac{h+1}{2}][P-\frac{l-h}{2}][P+\frac{l+h+2}{2}]}{\varphi_l(u-v)[P][P+h+1]}e^Q,\\
\widehat{\pi}_{l,w}(\hL^+_{+-}(u))&=&-S^-\frac{[u-v+\frac{h-1}{2}+P][\frac{l-h+2}{2}]}{\varphi_l(u-v)[P+h-1]}e^{-Q},\\
\widehat{\pi}_{l,w}(\hL^+_{-+}(u))&=&S^+\frac{[u-v-\frac{h+1}{2}-P][\frac{l+h+2}{2}]}{\varphi_l(u-v)[P]}e^Q,\\
\widehat{\pi}_{l,w}(\hL^+_{--}(u))&=&-\frac{[u-v-\frac{h-1}{2}]}{\varphi_l(u-v)}e^{-Q},
\en
where we set $z=q^{2u}$, $w=q^{2v}$, and  
\be
\varphi_l(u)&=&-z^{-{l}/{2r}}\rho_{1l}^+(z,p)^{-1}[u+\frac{l+1}{2}],\\
\rho_{kl}^+(z,p)&=&q^{{kl}/{2}}\frac{\{pq^{k-l+2}z\}\{pq^{-k+l+2}z\}}{\{pq^{k+l+2}z\}\{pq^{-k-l+2}z\}}\frac{\{q^{k+l+2}/z\}\{q^{-k-l+2}/z\}}{\{q^{k-l+2}/z\}\{q^{-k+l+2}/z\}}.
\en
\end{thm}

The following Proposition indicates a consistency of our construction of $\widehat{\pi}_{l,w}$ 
 and the fusion construction of the 
dynamical $R$ matrices (=face type Boltzmann weights).  
\begin{prop}\lb{LandR}
Let us define the matrix elements of $\widehat{\pi}_{l,w}(\hL^+_{\vep_1\vep_2}(u))$ by
\be
\widehat{\pi}_{l,w}(\hL^+_{\vep_1\vep_2}(u))v^l_m&=&
\sum_{m'=0}^l (\hL^+_{\vep_1\vep_2}(u))_{\mu_{m'}\mu_m}v^l_{m'},
\en
where $\mu_m=l-2m$. Then we have
\be
(\hL^+_{\vep_1\vep_2}(u))_{\mu_{m'}\mu_{m}}=
R^+_{1l}(u-v,P)_{\vep_1 \mu_{m'} }^{\vep_2 \mu_m}.
\en
Here $R^+_{1l}(u-v,P)$ is the $R$ matrix from (C.17) in  \cite{JKOS2}. 
The  case $l=1$, $R^+_{11}(u-v,P)$ coincides with the image 
$(\pi_{1,z}\otimes \pi_{1,w})$ of the universal $R$ matrix $\cR^+(\la)$\cite{JKOS} given in \eqref{Rmat}.
The case $l>1$, $R^+_{1l}(u-v,P)$ coincides with the $R$ matrix 
obtained by fusing $R^+_{11}(u-v,P)$ $l$-times. In particular the matrix element $R^+_{1l}(u-v,P)_{\vep \mu}^{\vep'\mu'}$ is gauge equivalent to the fusion face weight $W_{l1}(P+\vep',P+\vep'+\mu',P+\mu,P|u-v)$ from (4) in \cite{DJMO}. 
\end{prop}

\subsection{Infinite Dimensional Representation}

Let $V(\la_l)$ be the level-$k$ $(c=k)$ irreducible highest weight 
$\FF[U_q(\slth)]$-module
 of highest weight $\la_l=(k-l)\Lambda_0+l\Lambda_1\ (0\leq l\leq k)$. 
Here $\Lambda_i\ (i=0,1)$ denote the fundamental weights of $\slth$. 
We regard $\widehat{V}(\la)=\bigoplus_{m\in \Z}V(\la)\otimes \C e^{-mQ}$ as
the $U_{q,p}(\slth)$-module \cite{JKOS2}. 

We realize $\widehat{V}(\la_l)$ by using the Drinfeld generators $a_n\ (n\in \Z_{\not=0})$ and the $q$-deformed $\Z_k$-parafermion algebra\cite{Konno,JKOS2,KKW}.  Let us define 
$\al_n~(n \in {\mathbb Z}_{\neq 0})$ by 
\begin{eqnarray*}
\al_n=\left\{
\begin{array}{cc}a_n
& {\rm for }\ n>0\\
\frac{[rn]_q}{[r^* n]_q}q^{k|n|}a_n&{\rm for }\  n<0
\end{array}
\right.
\end{eqnarray*}
with $r^*=r-k$. Then we have 
\begin{eqnarray*}
~[\al_m,\al_n]=\frac{[2m]_q[km]_q}{m}
\frac{[rm]_q}{[r^*m]_q}\delta_{m+n,0}.
\end{eqnarray*}
The $q$-deformed $\Z_k$-parafermion algebra is an associative 
algebra over $\C$ generated by 
$\Psi_{+,\frac{\mu}{k}-n},\ \Psi_{-,\frac{\mu}{k}-n} \ 
(\mu, n\in \Z) $. Consider the generating functions (parafermion fields)
\be
&&\Psi(z)\equiv\Psi^+(z)=\sum_{n\in \Z}\Psi_{+,\frac{\mu}{k}-n} z^{-{\mu}/{k}+n-1},\\
&&\Psi^\dagger(z)\equiv\Psi^-(z)=\sum_{n\in \Z}\Psi_{-,\frac{\mu}{k}-n} 
z^{{\mu}/{k}+n-1}
\en
defined on a weight vector $v$ satisfying $q^hv=q^\mu v$. 
The parafermion fields $\Psi(z)$ and $\Psi^\dagger(z)$ satisfy
\be
&&\left(\frac{z}{w}\right)^{{2}/{k}}\frac{(x^{-2}w/z;x^{2k})_{\infty}}
{(x^{2+2k}w/z;x^{2k})_{\infty}}\Psi^{\pm}(z)\Psi^{\pm}(w)=
\left(\frac{w}{z}\right)^{{2}/{k}}\frac{(x^{-2}z/w;x^{2k})_{\infty}}
{(x^{2+2k}z/w;x^{2k})_{\infty}}\Psi^{\pm}(w)\Psi^{\pm}(z),\\
&&\left(\frac{z}{w}\right)^{-{2}/{k}}\frac{(x^{2+k}w/z;x^{2k})_{\infty}}
{(x^{-2+k}w/z;x^{2k})_{\infty}}\Psi^{\pm}(z)\Psi^{\mp}(w)-
\left(\frac{w}{z}\right)^{-{2}/{k}}\frac{(x^{2+k}z/w;x^{2k})_{\infty}}
{(x^{-2+k}z/w;x^{2k})_{\infty}}\Psi^{\mp}(w)\Psi^{\pm}(z)\nn\\
&&\qquad\qquad\qquad =\frac{1}{x-x^{-1}}\left(\delta\left(x^k\frac{w}{z}\right)-\delta\left(x^{-k}\frac{w}{z}\right)\right).
\en

\begin{thm}\cite{KKW}~~
By using the irreducible $q$-$\Z_k$ parafermion module $\H^{PF}_{l,M}$, the 
level-$k$ irreducible highest weight $U_{q,p}(\slth)$-module 
$\widehat{V}(\la_l)$ is realized as follows.
\be
&&\widehat{V}(\la_l)=\bigoplus_{m\in \Z}\bigoplus_{n\in\Z}
\bigoplus_{M\equiv 0\  {\rm mod}\ 2k\atop (M\equiv l\ {\rm mod}\ 2)}^{2k-1}\ \widehat{V}(\la_l)_{M+2kn+m},\\
&&\widehat{V}(\la_l)_{M+2kn+m}=\FF[\al_{-m}\ (m\in \Z_{>0})]\otimes 
\H^{PF}_{l,M}\otimes 
\C e^{{(M+2kn)}{\alpha}/{2}}\otimes \C e^{-mQ}. 
\en
The action of the elliptic 
 currents on $\widehat{V}(\la_l)$ are given by  
\begin{eqnarray*}
K(u)&\mapsto&:\exp\left(-
\sum_{m \neq 0}\frac{[m]_q}{[2m]_q[rm]_q}\al_{-m}z^m
\right):e^{Q}z^{-{k}({2}{P}-1)/{4rr^*}+{h}/{2r}},\\
E(u)&\mapsto&\Psi(z) :\exp\left(
-\sum_{m\neq 0}\frac{1}{[km]_q}\al_m z^{-m}
\right):e^{2Q+\alpha_1}
z^{{(h+1)}/{2}-{(P-1)}/{r^*}},\\
F(u)&\mapsto&\Psi(z)^\dagger:\exp\left(
\sum_{m\neq 0}\frac{[r^*m]_q}{[km]_q[rm]_q}\al_m z^{-m}
\right):e^{-\alpha_1}z^{-{(h-1)}/{2}
+{(P+h-1)}/{r}}.
\end{eqnarray*}
\end{thm}

Let $(\widehat{\pi}_V, V), (\widehat{\pi}_W, W)$ be two 
dynamical representations of $U_{q,p}$. We define the tensor product $V\tot W$ by
\be
&&V\tot W=\bigoplus_{\al\in \bar{\hh}^*}(V\tot W)_{\al},\quad 
(V \tot W)_{\al}=\bigoplus_{\beta\in \bar{\hh}^*}V_{\beta}\otimes_{M_{\h^*}}W_{\al-\beta}, 
\en
where
$\otimes_{M_{\h^*}}$ denotes the usual tensor product modulo the relation
\bea
&&f^*(P)v\otimes w=v\otimes f(P+h)w.\lb{VtotW}
\ena
Then,
 $(\widehat{\pi}_V\tot \widehat{\pi}_W)\circ \Delta : 
U_{q,p}\to \cD_{H,V}\tot \cD_{H,W}$ is a dynamical representation 
of $U_{q,p}$ on $V\tot W$.

\section{Vertex Operators}

By using the $\h$-Hopf algebroid structure, we define  
 the type I and II vertex operators of $U_{q,p}(\slth)$ 
as intertwiners of  $U_{q,p}(\slth)$-modules.  
Investigating their intertwining relations, we show that they 
 coincide with those obtained in \cite{JKOS2} by using the 
quasi-Hopf algebra structure of $\Bqla(\slth)$ and the isomorphism 
$U_{q,p}(\slth)\cong \Bqla(\slth)\otimes_{\C} \C[\bH^*]$. 

\begin{dfn}
The type I and II vertex operators of spin $n/2$ are the intertwiners of 
$U_{q,p}$-modules of the form
\be
\hPhi(u)\ &:&\ \widehat{V}(\la) \to  V^{(n)}_{z} \widetilde{\otimes}  \widehat{V}(\nu),\\
\hPsi^*(u)\ &:&\   \widehat{V}(\la)\widetilde{\otimes} V^{(n)}_{z} 
\to  \widehat{V}(\nu), 
\en
where $z=q^{2u}$, and $\widehat{V}(\la)$ and $\widehat{V}(\nu)$ denote the level-$k$ highest weight $U_{q,p}$-modules of  highest weights $\la, \nu$, respectively. 
They satisfy the
 intertwining relations with respect 
to the comultiplication $\Delta$ in Definition \ref{Hopfalgebroid}. 
\bea
&&\Delta(x)\hPhi(u)=\hPhi(u)x\qquad \forall x\in U_{q,p},\lb{typeIrel}\\
&&x\hPsi^*(u)=\hPsi^*(u)\Delta(x)\qquad \forall x\in U_{q,p}.\lb{typeIIrel}
\ena
\end{dfn}
The physically interesting cases are $n=k, \la=\la_l, \nu=\la_{k-l}$ for the type I and $n=1, \la=\la_l, \nu=\la_{l\pm1}$ for the type II. See for example \cite{KKW}. 

Let us define the components of the vertex operators as follows.
\bea
&&\hPhi(v-\frac{1}{2})=\sum_{m=0}^nv^{n}_m\widetilde{\otimes}\Phi_m(v),\lb{compoI}\\
&&\hPsi^*(v-\frac{c+1}{2})(\ \cdot\ \widetilde{\otimes} v_m^n)=\Psi^*_m(v).
\lb{compoII}
\ena 

\begin{thm}
The vertex operators satisfy the following linear equations. 
\bea
&&\hPhi(u)\hL^+(v)=R^{+(12)}_{1n}(v-u,P+h)\hL^+(v)\hPhi(u),\lb{typeIeq}\\
&&\hL^+(v)\hPsi^*(u)=\hPsi^*(u)\hL^+(v)R^{+*(13)}_{1n}(v-u,P-h^{(1)}-h^{(3)}).
\lb{typeIIeq}
\ena
The relation \eqref{typeIeq} 
should be understood  on $V^{(1)}_w\tot \widehat{V}(\la)$, whereas \eqref{typeIIeq} on 
$V^{(1)}_w\tot \widehat{V}(\la)\tot V^{(n)}_z$. 
\end{thm}

\noi
{\it Proof.} \ 
Applying $\Delta$ in Definition \ref{Hopfalgebroid}  
and noting Proposition \ref{LandR}, we obtain from \eqref{typeIrel}
\be
\hPhi(u)\hL^+_{\vep_1\vep_2}(v)
&=&\Delta(\hL^+_{\vep_1\vep_2}(v))\hPhi(u)\nn\\
&=&\sum_{m=0}^n\sum_{\vep}\hL^+_{\vep_1\vep}(v)v^n_m
\tot \hL^+_{\vep\vep_2}(v)
\Phi_m(u)\nn\\
&=&\sum_{m=0}^n\sum_{\vep}\sum_{m'=0}^n
R^+_{1n}(v-u,P)_{\vep_1 \mu_{m'}}^{\vep \mu_m}
v^n_{m'}\tot \hL^+_{\vep\vep_2}(v)
\Phi_m(u)\nn\\
&=&\sum_{m'=0}^nv^n_{m'}\tot \sum_{m=0}^n\sum_{\vep}
R^+_{1n}(v-u,P+h)_{\vep_1 \mu_{m'}}^{\vep \mu_m}\hL^+_{\vep\vep_2}(v)
\Phi_m(u),\lb{typeI}
\en
where $\mu_m=n-2m$ etc.
In the last equality we used \eqref{VtotW}. 
Similarly for the type II, from \eqref{typeIIrel}
we obtain  
\be
\hL^+_{\vep_1\vep_2}(u)\Psi^*_m(v+\frac{1}{2})
&=&\hPsi^*(v+\frac{1}{2})\Delta(\hL^+_{\vep_1\vep_2}(u))
(\ \cdot\ \widetilde{\otimes} v_m^n)\nn\\
&=&\sum_{\vep}
\sum_{m'}\hPsi^*(v+\frac{1}{2})\left(\hL^+_{\vep_1\vep}(u)\tot 
R^{+}_{1n}(u-v,P)_{\vep\ \mu_{m'}}^{\vep_2 \mu_{m}}v^n_{m'}\right)\nn\\
&=&\sum_{\vep}
\sum_{m'}\hPsi^*(v+\frac{1}{2})\left(R^{+*}_{1n}\left(u-v,P-\mu_{m'}\right)_{\vep\ \mu_{m'}}^{\vep_2 \mu_m}\hL^+_{\vep_1\vep}(u)\tot v^n_{m'}\right)\nn\\
&=&\sum_{\vep}
\Psi^*_{m'}(v+\frac{1}{2})
R^{+*}_{1n}\left(u-v,P-\mu_{m'}\right)_{\vep\ \mu_{m'}}^{\vep_2 \mu_m}
\hL^+_{\vep_1\vep}(u)\nn\\
&=&\sum_{\vep}
\Psi^*_{m'}(v+\frac{1}{2})\hL^+_{\vep_1\vep}(u)R^{+*}_{1n}\left(u-v,P-\mu_{m'}-\vep\right)_{\vep\ \mu_{m'}}^{\vep_2 \mu_m}.\lb{typeII}
\en
Here in the third equality, we used the relation \eqref{VtotW}. 
Note also $\vep+\mu_{m'}=\vep_2+\mu_m$.
\qed
 
\noi
The equations \eqref{typeIeq} and \eqref{typeIIeq} coincide with 
(5.3) and (5.4) in \cite{JKOS2}, respectively. 
Note that the comultiplication used in \cite{JKOS2} corresponds to
 the opposite one of $\Delta$ here. Under certain analyticity 
conditions, these equations 
determine the vertex operators uniquely up to normalization.

\subsection*{Acknowledgments}

The author would like to thank Michio Jimbo, Anatol Kirillov, Atsushi Nakayashiki,  Masatoshi Noumi and 
Hjalmar Rosengren for stimulating discussions and valuable suggestions. 
This work is supported by the Grant-in-Aid for Scientific Research (C)19540033, JSPS
Japan.


\end{document}